\documentclass[12pt]{article}
\usepackage[active]{srcltx}
\usepackage{amssymb, amsmath}
\begin{document}
\begin{center}
{\large\bf Phaseless inverse problems with interference waves}
\end{center}

\begin{center}
  V.~G.~ROMANOV\footnotemark[1] and  M.~YAMAMOTO \footnotemark[2]
\end{center}
\footnotetext[1] {Sobolev Institute of Mathematics,
              Siberian Division of Russian Academy of Sciences,
             Acad. Koptyug prospekt 4, 630090 Novosibirsk, Russia;
             e-mail: romanov@math.nsc.ru}
\footnotetext[2] {Department of Mathematical Sciences,
        The University of Tokyo,
        3-8-1 Komaba, Meguro, Tokyo, 153 Japan;
        e-mail: myama@ms.u-tokyo.ac.jp\\
Research Center of Nonlinear Problems of Mathematical Physics,
Peoples' Friendship University of Russia, Moscow, Russia
 }

\bigskip
{\bf Abstract.} {\footnotesize
We consider two phaseless inverse problems for elliptic equation.
The statements of these problems differ from have considered.  Namely,
instead of given information about modulus of scattering waves, we consider
the information related to modulus of full fields, which consist of sums of
incident and scattering fields. These full fields are the interference fields
generated by point sources.
We introduce a set of auxiliary point sources
for solving the inverse problems
and demonstrate that the corresponding data allow us to solve the inverse
problems in a way similar to the case of measurements of scattering waves.}



\bigskip
\section{Introduction}
\setcounter{equation}{0}

The phaseless inverse problem was stated first in the book \cite{ChSa} by
Chadan and Sabatier "Inverse Problems in Quantum Scattering Theory".
The authors noted that a phase of a solution of the Schr\"{o}dinger
equation cannot be measured for the large frequencies (energies).
Therefore they suggested to study inverse problems when only moduli of
the fields are measured.
The first results for the phaseless inverse problem related to recovering
a potential in the Schr\"{o}dinger equation were obtained in
the papers \cite{Klib1}-\cite{Klib3}, \cite{KR1}-\cite{KR3} by Klibanov and
Romanov and \cite{Nov1}-\cite{Nov3} by Novikov.
Then in the papers \cite{Klib4, KR4, KR5}  the phaseless problems were
considered for a generalized Helmholtz operator $\Delta +k^2n^2(x) $ with
incident point sources or incident plane waves.
In the paper \cite{KR6}, a procedure for extracting Riemannian distances $\tau(x,y)$ from given data of  the phaseless inverse problem was developed
for the Helmholtz operator.  Recently the inverse phaseless problem of recovering the permittivity coefficient in the Maxwell equations was considered in \cite{R3} for a point incident source and in \cite{R4} for an incident plane wave.

In the papers \cite{KR1} and \cite{KR5} by Klibanov and Romanov,
phaseless inverse problem for elliptic equations were studied.
The main features of the formulation of these problems consist in
the following: it is assumed that a coefficient of an elliptic equation
with a frequency parameter $k$ is known outside of some compact domain
$\Omega_0\subset\mathbb{R}^3$ with smooth boundary $S_0=\partial\Omega_0$
and it should be recovered in $\Omega_0$ from given modulus on $S_0$
of scattering fields generated by point sources $y\in S_0$ which
run over the whole $S_0$.  The given scattering field is then a function of
variable $x, y \in S_0$ and the frequency $k$. Then the asymptotic expansion
of the solutions to direct problems with a point source is studied
as $k\to\infty$.
As a result, the given information allows us to reduce the inverse problems
under consideration to well-known problems, that is,
the tomography problem for recovering a potential $q(x)$ in the
Schr\"{o}dinger equation (e.g., \cite{KR1}) or the inverse kinematic problem
for recovering a refractive index $n(x)$ in a generalized Helmholtz equation
(\cite{KR5}).
Related to this, we pose the question whether or not the modulus
not of the scattering field but a full field which is the the sum of
an incident field in the homogeneous medium and the scattering field in
non-homogeneities.
As it is seen from the above mentioned papers, given
scattering field simplifies study of the  phaseless inverse problem, but
the inverse problem with data of the full field is much more complicated.
Related results for an inverse problem with the full field were obtained
in a linearized approximation only (\cite{KR5}).

In this paper, we present a new approach for studying  phaseless inverse problems when the information in these problems is given for a full field.
In this approach, we suppose that
a potential or a refractive index is unknown only in the ball
$\Omega_0=\{x\in \mathbb{R}^3|\, |x|<R_0\}$ with $R_0>0$.
Let $\Omega=\{\mathbb{R}^3|\, |x|\le R\}$ be the ball with the radius
$R>R_0$.  We set $S=\{x\in \mathbb{R}^3|\, |x|=R\}$.
We consider an interference of two waves:
a wave produced by a point source located at $y\in S$ and
a wave located at the auxiliary point $z=z(y)\in\mathbb{R}^3\setminus
\Omega$.
The modulus of the field obtained by the interference of these waves is
measured at a receivers
$x\in S_+(z) := \{x\in S|\, x\cdot (z-x)>0\}$, and we regard the modulus
as a function in $x$, $y\in S$ and sufficiently large $k$.
Here we can consider $S_+(z)$ as the illuminated part of $S$ by the light
source placed at the point $z$.
Moreover, in our approach, for every $y\in S$ we need 3 different
sources $z^{(j)}=z^{(j)}(y)$, $j=1,2,3$.
It is assumed that these sources are on the tangent plane of $S$ at the point
$-y\in S$, lie enough far from the point  $-y$, and form a right
triangle on this tangent plane.  Hence the observation data are a function
of $y\in S$, $x\in S_+(z^{(j)})$ and $k\ge k_0>0$, where $k_0$ is a positive
number.  More detailed description is given in the next section.

The structure of this paper is as follows. In section 2,
we formulate our observation system which we use for consideration of two
phaseless inverse problems.
In section 3, we consider the Schr\"{o}dinger equation and the problem
of recovering a potential from phaseless data related to the full field.
Here we demonstrate that by given data, we can reduce the inverse problem
under consideration to a tomography problem.
In section 4,  we consider the problem of finding a refractive index $n(x)$
in a generalized Helmholtz equation from given modulus for 3 interference
fields and show that this problem is reduced to an inverse kinematic problem.
Thus it turns out that the phaseless data in our paper
allow us to reduce the inverse problems exactly to the same problems as in
\cite{KR1, KR5}.

\bigskip
\section{Formulation of the inverse problems: auxiliary point sources
for refractive waves}
\setcounter{equation}{0}
Here for the formulation of our inverse problems, we define an
observation system.
We shall consider differential equations with a source like
$\delta(x-y)+\delta(x-z)$, which is a superposition of
interference waves produced by two point sources $y, z \in \mathbb{R}$.

Let $\Omega_0 := \{x\in\mathbb{R}^3|\, |x|<R_0\}$ and
$\Omega=\{x\in\mathbb{R}^3|\, |x|<R\}$ with $R>R_0$, and $S =
\partial\Omega$.

Assume that $y$ is an arbitrary point of $S$ and for $y$, and we choose
a point $z=z(y)\in \mathbb{R}^3\setminus \Omega$ in the following way.
That is, consider the point $-y$. This point obviously belongs to $S$.
By $\Sigma(-y)$ we denote the tangent plane on $S$ at $-y$.
Let $z(y) \in \Sigma(-y)$ be a point satisfying
$\vert z(y) - (-y)\vert = \ell$ with
$\ell\ge \ell^*$, where  $\ell^*$ will be estimate below.
We set $S_+(z) := \{ x \in S\vert\thinspace x\cdot (z-x) > 0\}$, which can be
regarded as the part of $S$ illuminated by the light source placed at $z$.
Denote by $S^*(y)=\{x\in S|\, -y\cdot(x-y)>|x-y|\sqrt{R^2-R^2_0}\}$.
It means that $S^*(y)$ is the shadow part of $S$ for the light source placed
at $y$ when the light meets the non-transparent domain $\Omega_0$ on its way.
Later we need to construct for every
$y\in S$  some refractive waves with point
sources at $z^{(j)}(y)$, $j=1,2, \ldots,m$, such that the union  of $S_+(z^{(j)})$, $j=1,2, \ldots,m$, contains $S^*(y)$.
A simple analysis shows that it is impossible to do for arbitrary $R$ and $R_0$ if $m=1$ or $m=2$, but
it is possible if $m=3$ (see however Remark 1 below).  There are many possibilities in choosing such
$z^{(j)}(y)$, $j=1,2, 3$. Among them, we choose $z^{(j)}(y)$, $j=1,2, 3$
such that $z^{(j)}\in\Sigma(-y)$,
$|z^{(1)}+y|=|z^{(2)}+y|=|z^{(3)}+y|=\ell$
and the points $z^{(j)}(y)$, $j=1,2, 3$, are the vertices of a right triangle.
Then

 \bigskip
{\bf Lemma 1.} {\it Let   $z^{(j)}\in\Sigma(-y)$, $j=1,2, 3$, $|z^{(1)}+y|
=|z^{(2)}+y|=|z^{(3)}+y|=\ell$, and the points $z^{(j)}(y)$, $j=1,2, 3$ be
the vertices of a right triangle and $\ell\ge \ell^*$, where}
\begin{eqnarray}\label{4}
\ell^*=\frac{2R_0R}{\sqrt{R^2-R_0^2}}.
\end{eqnarray}
{\it Then  $S^*(y)\subset (S_+(z^{(1)})\cup S_+(z^{(2)})\cup S_+(z^{(3)}))$.}

 \bigskip
Proof.
By the invariance of $z^{(j)}$, $j=1,2,3$, it suffices to consider
the case where $y=(0,0,-R)$. Then $-y=(0,0,R)$.  On the plane $\Sigma(-y)$
we consider the points $z^{(1)}=(\ell/2,\ell\sqrt{3}/2, R)$,
$z^{(2)}=(-\ell/2,\ell\sqrt{3}/2, R)$, $z^{(3)}=(0,-\ell, R)$ and
demonstrate that the relation $S^*(y)\subset (S_+(z^{(1)})\cup
S_+(z^{(2)})\cup S_+(z^{(3)}))$ holds if $\ell\ge \ell^*$, where $\ell^*$
is defined by (\ref{4}).
Let $z=(\ell\cos\varphi, \ell\sin\varphi, R)$ with $\varphi\in(0,\pi)$,
and let $C(z) = \{ \xi \in S\vert \thinspace \vert \xi-z\vert =\ell\}$.

Then $\xi\in C(z)$, $\xi=(\xi_1,\xi_2,\xi_3)$ if and only if
$\xi=z+\ell\beta$ where $\beta=(\beta_1,\beta_2,\beta_3)$ is a unit
vector satisfying the conditions $z\cdot\beta=-\ell$ and $\beta_3\in (-1, 0]$.
The first of these conditions means that $|\xi|^2=R^2$, and
the second one  means that the component $\beta_3$ should be non positive.
Since $\varphi\in(0,\pi)$, the circumference $C(z)$
 has two points of intersections with the plane $\xi_1=0$.
One of these points is $-y=(0,0, R)$. Let us find the other intersection point. Note that for both intersection points the equality
$0=\ell\cos\varphi+\ell\beta_1$
holds.
Hence, $\beta_1=-\cos\varphi$. Then $\beta_3=-\sqrt{\sin^2\varphi-\beta^2_2}$. One can calculate $\beta_2$ using
 the relation $z\cdot\beta=-\ell$. Taking into account that
$z=(\ell\cos\varphi,\ell\sin\varphi,R)$,
 $\beta_1=-\cos\varphi$ and $\beta_3=-\sqrt{\sin^2\varphi-\beta^2_2}$, we
obtain the equality
$$
-\ell\cos^2\varphi+\beta_2 \ell\sin\varphi-R\sqrt{\sin^2\varphi-\beta^2_2}
= -\ell.
$$
From this, we have the quadratic equation in $\beta_2$:
$$
\beta_2^2 (R^2+\ell^2\sin^2\varphi)+2\beta_2 \ell^2\sin^3\varphi
+ \ell^2\sin^4\varphi-R^2\sin^2\varphi=0.
$$
Solving this equation, we obtain
 $$
(\beta_2)_{\pm}=\sin\varphi\frac{-\ell^2\sin^2\varphi\pm R^2 }
{R^2+\ell^2\sin^2\varphi}.
 $$
Here $(\beta_2)_{-}=-\sin\varphi$ corresponds to
the intersection point $-y$, while
$$
(\beta_2)_{+}=\sin\varphi\frac{R^2 - \ell^2\sin^2\varphi}
{R^2+\ell^2\sin^2\varphi}
 $$
corresponds to the second intersection point.
This point has the coordinates
\begin{equation}\label{5}
 \xi_1=0,\> \xi_2=\ell\sin\varphi+\ell(\beta_2)_{-},
\>\xi_3=R-\ell\sqrt{\sin^2\varphi-(\beta_2)_{-}^2}.
 \end{equation}
Hence,
\begin{equation}\label{6}
\xi_1=0, \quad \xi_2= \frac{2\ell R^2\sin\varphi}{R^2+\ell^2\sin^2\varphi},
\quad \xi_3=R-\frac{2R\ell^2\sin^2\varphi}{R^2+\ell^2\sin^2\varphi}.
\end{equation}
Obviously that $\xi_3\in (-R,R)$ for $\varphi\in (0,\pi)$.

Consider now two points $z^{(1)}=(\ell\cos(\pi/6), \ell\sin(\pi/6),R)$,
$z^{(2)}=(\ell\cos(5\pi/6), \ell\sin(5\pi/6),R)$.
These points lie symmetrically with respect to
the plane $\xi_1=0$. Therefore they have the same intersection points
 with the plane $\xi_1=0$, and  the coordinates of
$\xi^*\neq -y$ are determined by the formulae  (\ref{6}) with
$\varphi=\pi/6$, that is,
\begin{equation}\label{7}
 \xi_1^*=0, \quad
 \xi_2^*= \frac{2 \ell R^2}{4R^2+\ell^2},\quad
\xi_3^*=R-\frac{2R\ell^2}{4R^2+\ell^2}.
\end{equation}
Denote by $S_{12}(\ell)$ the piece of $S$ bounded
by the meridional semi-planes $\varphi=\pi/6$,
 $\varphi=5\pi/6$ and by the plane $\xi_3=\xi_3^*$. Obviously $S_{12}(\ell)
\subset (S_+(z^{(1)})\cup S_+(z^{(2)}))$.
Now we consider 3 points: $z^{(1)}=(\ell\cos(\pi/6), \ell\sin(\pi/6), R)$,
$z^{(2)}=(\ell\cos(5\pi/6), \ell\sin(5\pi/6),R)$  and
  $z^{(3)}=(\ell\cos(3\pi/2),\ell\sin(3\pi/2), R)$. All these points are
in a symmetrical position  with respect to one to other.
Therefore if we denote by $S_{13}(\ell)$ the part of $S$ bounded by the
meridional semi-planes $\varphi=\pi/6$,
 $\varphi=3\pi/2$ and by the plane $\xi_3=\xi_3^*$,
then we find that $S_{13}(\ell)\subset (S_+(z^{(1)})\cup S_+(z^{(3)}))$.

Similarly, if $S_{23}(\ell)$ denotes the part of $S$ bounded by the meridional
semi-planes $\varphi=5\pi/6$, $\varphi=3\pi/2$ and by the plane
$\xi_3=\xi_3^*$, then $S_{23}(\ell)\subset (S_+(z^{(1)})\cup S_+(z^{(3)}))$.
Hence $S(\ell) :=  (S_{12}(\ell)\cup S_{13}(\ell)\cup S_{23}(\ell))
\subset (S_+(z^{(1)})\cup S_+(z^{(2)})\cup S_+(z^{(3)}))$.
 On the other hand, $S(\ell)=\{\xi\in\mathbb{R}^3|\, |\xi|= R,
\xi_3\ge \xi_3^*\}$.  Note that $\xi_3^*\to -R$ as $\ell\to\infty$.
 Choose $\ell$ such that $S^*(y)\subset S(\ell)$. For this we first need
find  $h=\inf \xi_3$ for all $\xi\in S^*(y)$.
By $\gamma$ we denote the angle between the $\xi_3$-axis
and a straight line which passes $y$ and is tangent on $\Omega_0$.
Then $\sin\gamma=  R_0/R$. The length of the piece of this straight line
included in $\Omega$ is $2R\cos\gamma$ and $h=2R\cos^2\gamma-R$.
Making simple calculations we find
$$
h=2R\left(1-\frac{R_0^2}{R^2}\right)-R.
$$
The inclusion $S^*(y)\subset S(\ell)$ holds if $h\ge \xi_3^*$, that is,
$$
2R\left(1-\frac{R_0^2}{R^2}\right)-R \ge R-\frac{2R\ell^2}{4R^2+\ell^2}.
$$
The latter inequality holds if $\ell\ge \ell^*$ where $\ell^*$ is given
by (\ref{4}).
\hfill $\Box$

The points $z^{(j)}(y)$ satisfying conditions of
Lemma 1 specify the auxiliary  points sources for measuring phaseless data
related to the source $y$.

 \bigskip
{\bf Remark 1.} If $R>\sqrt{2} R_0$ then instead of 3 points $z^{j}, j=1,2,3,$ one can use only one point $z=-\lambda y$ with $\lambda\ge R^2/(R^2-\sqrt{2} R_0^2)$. In this case $S^*(y)\subset S_+(z)$.  Thus in this case one auxiliary point  $z=-\lambda y$ and the data on the full  interference field corresponding
 point sources placed at $y$ and $z$ play  the same role as the points $z^{j}, j=1,2,3,$   and data on 3 full  interference fields related to them.

\bigskip
\section{The phaseless inverse problem of determining a potential
for the Schr\"{o}dinger equation}
\setcounter{equation}{0}
Throughout this paper, let $i = \sqrt{-1}$ and
$\delta$ be the Dirac delta function.

We study here the phaseless inverse problem for  the  equation
\begin{eqnarray}\label{1}
-(\Delta +k^2-q(x))u  =\delta(x-y)+\delta(x-z), \>x\in\mathbb{R}^3.
\end{eqnarray}
Let a function $u=u(x,y,z,k)$ satisfy the equation (\ref{1})
and the radiation conditions
\begin{eqnarray}\label{2}
u=O(r^{-1}), \quad  \frac{\partial u}{\partial r}- iku=o(r^{-1}) \quad \text{as}\> r=|x| \to\infty.
\end{eqnarray}
Here  conditions (\ref{2}) is assumed to be valid uniformly
for all the directions $x/r$.

We assume that a potential $q(x)$ satisfies
\begin{eqnarray}\label{3}
q\in C^{4}(\mathbb{R}^3), \quad \> q(x)\ge0, \thinspace
x\in \mathbb{R},  \quad \text{supp}\, q \subset\Omega_0.
\end{eqnarray}
The solution of the problem (\ref{1}) and (\ref{2}) describes an interference
of waves produced by the sources at points $y$ and $z$.
We suppose that the modulus of the function $u(x,y,z,k)$ can be measured on the sets
$S_+(z^{(j)})$, $ j=1,2,3$,  introduced above.

\medskip
{\bf Problem 1.} Let $y\in S$ and $z^{(j)}=z^{(j)}(y)$, $j=1,2,3$ be
defined as above.
Then find $q(x)$ in $\Omega_0$ from
\begin{eqnarray}\label{8}
|u(x,y, z^{(j)}, k)|^2= f_j(x,y,k),\quad y\in S, \> x\in S_+(z^{(j)}), \> k\ge k_0,\> j=1,2,3,
\end{eqnarray}
where $u(x,y, z^{(j)}, k)$ is the solution to the problem (\ref{1})-(\ref{2})
with $z=z^{(j)}$ and $k_0>0$ is a fixed positive number.

\medskip
Our first main result is the uniqueness in the inverse problem.
\\
{\bf Theorem 1.} {\it Let a potential $q(x)$ satisfy (\ref{3}) and
the additional condition $|z^{(j)}(y)+y|=\ell>R\sqrt{8}$, $ j=1,2,3$.
Then $q$ is uniquely determined by the information (\ref{8}).}

\medskip
Proof. Represent the solution   $u(x,y,z,k)$ of the problem (\ref{1})-(\ref{2}) in the form
\begin{eqnarray}\label{9}
u(x,y, z, k)= v(x,y,k)+v(x,z,k),
\end{eqnarray}
where $v(x,y,k)$ satisfies the equation
\begin{eqnarray}\label{10}
-(\Delta +k^2-q(x))v  =\delta(x-y), \>x\in\mathbb{R}^3
\end{eqnarray}
and the radiation conditions.
Then $v(x,y,k)$  can be represented in  the form
\begin{eqnarray}\label{11}
v(x,y,k)=u_0(x,y,k)+v_{sc}(x,y,k),
\end{eqnarray}
where $v_0(x,y,k)$ is defined by the formula
\begin{eqnarray}\label{12}
v_0(x,y,k)=A_0(x,y)e^{ik|x-y|}, \quad A_0(x,y)=\frac{1}{4\pi |x-y|},
\end{eqnarray}
and is the fundamental solution of the Helmholtz operator $-\Delta-k^2$ with
the conditions (\ref{2}) and $v_{sc}(x,y,k)$ is the scattering field from the
point source located at $y$.

Consider the asymptotic behavior of $v_{sc}(x,y,k)$ as $k\to\infty$.
Theorem 1 in \cite{KR1} yields

\medskip
{\bf Lemma 2}. {\it Suppose that a potential $q(x)$ satisfies
conditions (\ref{3}).  Then for each pair of points
$x,y$ with $x\neq y$ the asymptotic behavior of the function
$v_{sc}(x,x^{0},k)$ is
\begin{equation}\label{13}
v_{sc}( x,y,k) =\frac{i\exp (ik|x-y|)}{8\pi |x-y|k}
\left[\int\limits_{L( x,y)} q( \xi ) d\sigma +O\Big(\frac{1}{k}\Big)\right] ,\>k\rightarrow \infty,
\end{equation}%
where $L( x,y)$ is the segment of the straight line connecting $x$ and $y$,
and $d\sigma$ is the arc length.}

\medskip
By this lemma, we can obtain the asymptotic behavior of the data
(\ref{8}):
\begin{eqnarray*}
 f_j(x,y,k)=\left|A_0(x,y)e^{ik|x-y|}+\frac{ie^{ik|x-y|}}{8\pi |x-y|k}
\int\limits_{L( x,y)} q( \xi ) d\sigma +A_0(x,z^{(j)}(y))e^{ik|x-z^{(j)}(y)|}\right|^2+O\Big(\frac{1}{k^2}\Big)
\\
=A^2_0(x,y)+A^2_0(x,z^{(j)}(y))+2A_0(x,y)A_0(x,z^{(j)}(y))\cos(k|x-z^{(j)}(y)|-k|x-y|)\\
\\
+\frac{A_0(x,z^{(j)}(y))}{4\pi |x-y|k}\int\limits_{L( x,y)} q( \xi ) d\sigma
\sin(k|x-z^{(j)}(y)|-k|x-y|)+O\Big(\frac{1}{k^2}\Big),
 \\
  y\in S, \> x\in S_+(z^{(j)}), \> k\ge k_0, \> j=1,2,3.
\end{eqnarray*}
For $x\in S_+(z^{(j)})$, we have
$$
\min_{x\in S_+(z^{(j)})}|x-z^{(j)}(y)|\ge \sqrt{\ell^2+R^2}-R, \quad
\max_{x\in S_+(z^{(j)})}|x-y|\le 2R.
$$
Therefore under the condition $\sqrt{\ell^2+R^2}-R>2R$, that is,
$\ell>R\sqrt{8}$,
the inequality $|x-z^{(j)}(y)|-|x-y|>0$ holds  for all $x\in S_+(z^{(j)})$.
Fix $y$ and $x\in S_+(z^{(j)})$ and take $k=k_m(x,y)=\pi/2+2\pi/(|x-z^{(j)}(y)|-|x-y|)$.  Then we obtain
 \begin{eqnarray*}
 f_j(x,y,k_m(x,y)) =A^2_0(x,y)+A^2_0(x,z^{(j)}(y))
 \qquad\qquad\qquad\qquad\qquad\qquad
 \\
 + \frac{A_0(x,z^{(j)}(y))}{4\pi |x-y|k_m(x,y)}\int\limits_{L( x,y)} q( \xi ) d\sigma+
+O\Big(\frac{1}{k_m^2(x,y)}\Big),
\\
  y\in S, \> x\in S_+(z^{(j)}), \> k_m(x,y)\ge k_0, \> j=1,2,3.
\end{eqnarray*}
Hence,
\begin{eqnarray}\label{14}
\int\limits_{L( x,y)} q( \xi ) d\sigma=g_j(x,y), \quad y\in S, \> x\in S_+(z^{(j)}), \> j=1,2,3,
\end{eqnarray}
where
\begin{eqnarray*}
g_j(x,y)=4\pi |x-y|\lim_{k_m(x,y)\to\infty}\frac{ [
 f_j(x,y,k_m(x,y)) -A^2_0(x,y)-A^2_0(x,z^{(j)}(y))] k_m(x,y)}{A_0(x,z^{(j)}(y))},
 \\
  y\in S, \> x\in S_+(z^{(j)}), \> j=1,2,3.
 \end{eqnarray*}
Because $S^*(y)\subset (S_+(z^{(1)})\cup S_+(z^{(2)})\cup S_+(z^{(3)}))$,
the equality yields that integrals over $q(x)$ along
 $L(x,y)$ are given for all $y\in S$ and $x\in S^*(y)$.
Hence we know the integrals for all straight lines that cross
 out the ball $\Omega_0$.
Hence we reach a usual tomography problem for recovering $q(x)$ inside
$\Omega_0$.
It proves Theorem 1.
 \hfill $\Box$

\medskip
We see that the given information related to full fields for reflected waves
allows us to reduce the phaseless inverse problem to the same tomography
problem as for data related to simple sources at $y$  and the scattering filed.

The same idea successively works in many others phaseless inverse problems.
Below we consider such an inverse problem for the Helmholtz equation.

\bigskip
\section{The phaseless inverse problem of determining a refractive
index for the Helmholtz equation}
\setcounter{equation}{0}
Following \cite{KR5}, we consider the phaseless inverse problem for the
Helmholtz equation.
Let $\Omega_0$, $\Omega $ and $S$ be as in the previous section.
Let $n(x)$, $x\in \mathbb{R}^{3}$ be a real valued function satisfying
\begin{equation}\label{15}
n\in C^{15}(\mathbb{R}^{3}),\quad 1\le n(x)\le n_1< \infty, \thinspace
x \in \mathbb{R}, \quad
\text{supp}(n-1)\subset\Omega_0,
\end{equation}
where $n_1$ is a given constant.
Let $u(x,y,z,k)$ satisfy
\begin{equation}\label{16}
-(\Delta+k^2n(x))u=\delta(x-y)+\delta(x-z), \quad x\in \mathbb{R}^3,
\end{equation}
and the radiation condition (\ref{2}).
Consider the following problem

\medskip
{\bf Problem 2.} Let $y\in S$ and $z^{(j)}=z^{(j)}(y)$, $j=1,2,3$, be
defined as above.
Find $n(x)$ in $\Omega_0$ from the data
\begin{eqnarray}\label{17}
|u(x,y, z^{(j)}, k)|^2= f_j(x,y,k),\quad y\in S,
\> x\in S_+(z^{(j)}), \> k\ge k_0,\> j=1,2,3,
\end{eqnarray}
where $u(x,y, z^{(j)}, k)$ is the solution to problem (\ref{16}) and
(\ref{2}) with $z=z^{(j)}$ and $k_0>0$ is a fixed positive number.

\medskip
Introduce the conformal Riemannian metric by the formula
$$
d\tau=n(x) |dx|, \quad  |dx|=\left(\sum_{i=1}^3 dx_i^2\right)^{1/2},
$$
where $d\tau$ is the element of length and by $\tau(x,y)$ we denote
the Riemannian distance between the points $x, y \in \mathbb{R}^3$.
As in \cite{KR5}, we shall pose the following assumption.

\medskip
{\bf Assumption}. {\it The Riemannian metric $d\tau = n(x)\vert dx\vert$
is simple, that is,
every two points $x,y\in \mathbb{R}^{3}$  can be connected by a single
geodesic line $\Gamma(x,y)$. }

\medskip

\medskip
{\bf Theorem 2.} {\it  Let conditions (\ref{15}) and Assumption
be fulfilled.  Additionally we assume that $z^{(j)}(y)$ satisfies
$ |z^{(j)}(y)+y|=\ell> R\sqrt{(1+2n_1)^2-1}$, $j=1,2,3$.
Then the data (\ref{17}) uniquely determines $n(x)$ inside $\Omega_0$. }

\medskip
{\bf Proof.}
We use the following lemma which is a corollary of Theorem 2
and formula (3.18) from \cite{KR5}

\
\medskip
{\bf Lemma 3.} {\it  Let conditions (\ref{15}) and Assumption  be fulfilled.
If $v=v(x,y,k)$ satisfies
\begin{equation}\label{18}
-(\Delta+k^2n(x))v=\delta(x-y), \quad x\in \mathbb{R}^3
\end{equation}
and the radiation condition
\begin{eqnarray}\label{19}
v=O(r^{-1}), \quad  \frac{\partial v}{\partial r}- ikv=o(r^{-1})
\quad \mbox{as $r=|x| \to\infty$}
\end{eqnarray}
for $x\neq y$, then we have the asymptotic formula
\begin{eqnarray}\label{20}
v(x,y,k)=A(x,y) e^{ik\tau(x,y)} +O\left(\frac{1}{k}\right) \quad
\> \text{as} \>\> k\to\infty,
\end{eqnarray}
where $A=A(x,y)>0$, $\in C^{13}((\mathbb{R}^3\setminus\{y\})\times S)$.
}

\medskip

It follows from Lemma 3 that the solution of the problem (\ref{16}) and
(\ref{2}) for $x\neq y$ posseses the asymptotic behavior of the form
\begin{equation}\label{21}
u(x,y,z,k)=A(x,y) e^{ik\tau(x,y)}+A(x,z) e^{ik\tau(x,z)}+O\left(\frac{1}{k}
\right) \quad \mbox{as $k\to\infty$}.
\end{equation}
Hence we conclude that
\begin{eqnarray}\label{22}
  f_j(x,y,k)=|u(x,y,z^{(j)}(y),k)|^2  =A^2(x,y) +A^2(x,z^{(j)}(y))
 \nonumber \\
  +2A(x,y)A(x,z^{(j)}(y))\cos(k\tau(x,z^{(j)}(y))-k\tau(x,y))
  +O\left(\frac{1}{k}\right)              \\
\> \text{as} \>\> k\to\infty, \quad y\in S, \quad x\in S_+(z^{(j)}(y)).
\nonumber
\end{eqnarray}
Note that $\tau(x,z^{(j)}(y))=|x-z^{(j)}(y)|$
for $y\in S$ and $x\in S_+(z^{(j)}(y))$. Moreover,
$$
\min\tau(x,z^{(j)}(y))\ge \sqrt{\ell^2+R^2}-R,
$$
where $\ell=|z^{(j)}(y)+y|$. On the other hand, for $y\in S$ and
$x\in S_+(z^{(j)}(y))$ we have $\tau(x,y)\le 2Rn_1$, where
$n_1$ is the bound given in (\ref{15}).

The latter equality holds because for fixed $x$ and $y$, the function
$\tau(x,y)$ minimizes integrals of $n(x)$ over arbitrary
smooth curves connecting $x$ and $y$.
In particular, $\tau(x,y)$ is less than the integral along the straight
line $L(x,y)$.
Therefore $\tau(x,y)\le n_1|x-y|\le 2n_1R$.
Since the condition in Theorem 2 implies $\ell>R\sqrt {(1+2n_1)^2-1}$,
we have
$$
\rho(x,y)=\tau(x,z^{(j)}(y))-\tau(x,y)\ge \sqrt{\ell^2+R^2}-R-2Rn_1>0.
$$
Fix $x$ and $y$ in (\ref{22}).
Then the left-hand side is an almost periodic function of $k$.
Hence, we can extract
the period of this function and find the difference
$\rho(x,y)=\tau(x,z^{(j)}(y))-\tau(x,y)$. The procedure of extracting the
period of almost periodic function is given in the paper \cite{KR6}.
Then one can calculate
$\tau(x,y)=\tau(x,z^{(j)}(y))-\rho(x,y)$ for all $y\in S$ and
$x\in S_+(z^{(j)}(y))$. Since
$S^*(y)\subset(S_+(z^{(1)}(y)\cup S_+(z^{(2)}(y)\cup S_+(z^{(3)}(y))$,
we find $\tau(x,y)$ for all $y\in S$ and $x\in S^*(y)$.
Note that $\tau(x,y)=|x-y|$ for all $y\in S$ and $x\in S\setminus S^*(y)$.
Therefore $\tau(x,y)$ is known for all
$(x,y)\in S\times S$.
Therefore we reach the well-known inverse kinematic problem:
find $n(x)$ in $\Omega_0$ from given $\tau(x,y)$ for all
$(x,y)\in S\times S$.

The multidimensional inverse kinematic problem  was studied  for first time
in a linear approximations in
\cite{LR, R1, R2}. The nonlinear problem was studied under Assumption
in the papers \cite{BG, Mukh, MR}, where uniqueness and stability theorems were shown. From the results given in these papers, the proof of Theorem 2 is
complete.

\hfill $\Box$

\medskip
\begin{center}
\textbf{Acknowledgments}
\end{center}
The work of V. G. Romanov was partially supported by the Russian Foundation for Basic Research grant
No. 17-01-00120.
The work of M. Yamamoto was supported partly by
Grant-in-Aid for Scientific Research (S) 15H05740 of
Japan Society for the Promotion of Science and
the Ministry of Education and Science of the Russian Federation (the
Agreement number No. 02 a03. 21. 0008).

\end{document}